\long\def\ignore#1{}
\font\largebf=cmbx10 scaled\magstep2
\font\tenmsa=msam10
\font\tenmsb=msbm10
\def\aff{\hbox{\rm aff}}
\def\all{\hbox{for all}}
\def\and{\hbox{and}}

\def\bigcupn{\bigcup\nolimits}
\def\caff{\hbox{\rm caff}}
\def\cite#1\endcite{[#1]}
\def\dom{\hbox{\rm dom}}

\def\dst{^{**}}
\def\eps{\varepsilon}

\def\exs{\hbox{there exists}}

\def\infn{\inf\nolimits}
\def\intr{{\rm int}}

\def\KST{K\"onig's sandwich theorem}
\def\lr{\Longrightarrow}
\def\minn{\min\nolimits}

\def\on{\hbox{on}}
\def\PC{{\cal PC}}
\def\PCLSC{{\cal PCLSC}}
\def\phi{\varphi}
\def\qed{\hfill\hbox{\tenmsa\char03}}
\def\qlr{\quad\lr\quad}
\def\quand{\quad\and\quad}
\def\r{\hbox{\tenmsb R}}
\def\rbar{\,]{-}\infty,\infty]}

\long\def\slant#1\endslant{{\sl#1}}
\def\st{\hbox{such that}}
\def\supn{\sup\nolimits}
\def\T{{\cal T}}
\def\ts{\textstyle}
\def\wh{\widehat}

\def\Zh{{\wh Z}}
\def\({\big(}
\def\){\big)}
\def\[{\big[}
\def\]{\big]}
\def\defSection#1{}
\def\defCorollary#1{}
\def\defDefinition#1{}
\def\defExample#1{}
\def\defLemma#1{}
\def\defNotation#1{}
\def\defProblem#1{}
\def\defRemark#1{}
\def\defTheorem#1{}
\def\locno#1{}
\def\meqno#1{\eqno(#1)}
\def\nmbr#1{}
\def\Proof{\medbreak\noindent{\bf Proof.}\enspace}
\def\Proo{{\bf Proof.}\enspace}
\def\Signoff{}
\def \INTsec{1}
\def \HBLsec{2}
\def \PRELIMdef{1}
\def \SUBdef{2}
\def \HBlem{3}
\def \ASANDthm{4}
\def \AFFSANDone{1}
\def \CHOICErem{5}
\def \STARdef{6}
\def \SUBLEVthm{7}
\def \ZXone{2}
\def \ZXtwo{3}
\def \ZXthree{4}
\def \ZXnew{5}
\def \ZXfour{6}
\def \ZXfive{7}
\def \ZXsix{8}
\def \ZXseven{9}
\def \SUBFUNDrem{8}
\def \TRIcor{9}
\def \TRIone{10}
\def \TRItwo{11}
\def \BDsec{3}
\def \BDTRIthm{10}
\def \BDZPXtwo{12}
\def \BDZPXfour{13}
\def \BDZPXthree{14}
\def \BDFENcor{11}
\def \BDFENtwo{15}
\def \SANDFENrem{12}
\def \BDQUADthm{13}
\def \BDQUADtwo{16}
\def \BDQUADthree{17}
\def \BDBIBIcor{14}
\def \BDBIBItwo{18}
\def \BDBIBIthree{19}
\def \BDSZcor{15}
\def \BDINDcor{16}
\def \BDINDone{20}
\def \BDINDtwo{21}
\def \BDINDthree{22}
\def \BDINDfour{23}
\def \BDINDSPECcor{17}
\def \FNORMsec{4}
\def \CAUCHYlem{18}
\def \FNORMlem{19}
\def \RSthm{20}
\def \RSone{24}
\def \RStwo{25}
\def \RSthree{26}
\def \RSfour{27}
\def \RSfive{28}
\def \RSsix{29}
\def \RScor{21}
\def \COMPLsec{5}
\def \TRIthm{22}
\def \ZPXtwo{30}
\def \ZPXthree{31}
\def \FENcor{23}
\def \FENtwo{32}
\def \QUADthm{24}
\def \QUADtwo{33}
\def \QUADthree{34}
\def \BIBIcor{25}
\def \BIBIone{35}
\def \SZthree{36}
\def \SZcor{26}
\def \EEAcor{27}
\def \EEBcor{28}
\def \INDcor{29}
\def \INDSPECcor{30}
\def \BANACH{1}
\def \KN{2}
\def \KONIG{3}
\def \KOTHE{4}
\def \RTRFENCHEL{5}
\def \RS{6}
\def \RUDIN{7}
\def \MSF{8}
\def \HBL{9}
\def \HBM{10}
\def \QUADARCHIV{11}
\def \QUADARCHIVTWO{12}
\def \SZNZ{13}
\def \ZIA{14}
\def \ZZOR{15}
\def \ZBOOK{16}
\magnification 1200
\headline{\ifnum\folio<2
{\hfil{\largebf The asymmetric sandwich theorem}\hfil}
\else\centerline{\bf The asymmetric sandwich theorem}\fi}
\medskip
\centerline{\bf S. Simons}
\medskip
\centerline{\sl Department of Mathematics, University of California}
\centerline{\sl Santa Barbara, CA 93106-3080}
\centerline{\sl simons@math.ucsb.edu}
\medskip
\centerline{\bf Abstract}
\medskip
We discuss the asymmetric sandwich theorem, a generalization of the Hahn--Banach theorem.   As applications, we derive various results on the existence of linear functionals that include bivariate, trivariate and quadrivariate generalizations of the Fenchel duality theorem.   Most of the results are about affine functions defined on convex subsets of vector spaces, rather than linear functions defined on vector spaces.   We consider both results that use a simple boundedness hypothesis (as in Rockafellar's version of the Fenchel duality theorem) and also results that use Baire's theorem (as in the Robinson--Attouch--Brezis version of the Fenchel duality theorem).   This paper also contains some new results about metrizable topological vector spaces that are not necessarily locally convex.
\defSection \INTsec
\bigbreak
\centerline{\bf\INTsec.\quad Introduction}
\medskip
\noindent
This paper is about the existence of linear functionals in various situations.   The main results of Section \HBLsec\ are the \slant asymmetric sandwich theorem\endslant\ of Theorem \ASANDthm\ and the \slant sublevel set theorem\endslant\ of Theorem \SUBLEVthm.   The asymmetric sandwich theorem is a straightforward extension of the Hahn--Banach theorem, and the sublevel set theorem is about the conjugate of a proper convex function defined on a CC space (a nonempty convex subset of some vector space).
\par
The final result in Section \HBLsec, Corollary \TRIcor, is a technical result about the conjugate of a convex function defined as the infimum of another convex function over a variable set.   We call this result a \slant trivariate existence theorem \endslant\ because it uses three spaces: two CC spaces and a locally convex space.   The best way of clarifying the interiority condition (\TRIone) that appears in Corollary \TRIcor\ is to consider the case of two generalizations to infinite dimensional spaces of the Fenchel duality theorem for two convex functions $f$ and $g$.   The first, due to Rockafellar, assumes a local boundedness condition for $g$ at some point where $f$ is finite, but does not assume any global lower semicontinuity conditions for either function.   Furthermore, the two functions are not treated symmetrically.   The second, due to Robinson and Attouch--Brezis, treats the two functions in a symmetric fashion.   Furthermore, two functions are assumed to be lower semicontinuous and the space complete. 
\par
Corollary \TRIcor\ leads to further trivariate existence theorems:  Theorem \BDTRIthm\ in Section \BDsec, and Theorem \TRIthm\ in Section \COMPLsec.   At this point we will discuss Theorem \BDTRIthm, since our remarks about Theorem \TRIthm\ are best postponed until after our consideration of Section \FNORMsec.   We give two consequences of Theorem \BDTRIthm:  Corollary \BDFENcor, an affine Fenchel duality theorem, and Theorem \BDQUADthm.  We call Theorem \BDQUADthm\ a \slant quadrivariate existence theorem \endslant\ because it uses four spaces: two CC spaces, a vector space and a locally convex space.   Corollary \BDFENcor\ extends the result of Rockafellar referred to above, and Theorem \BDQUADthm\ shows how we can compute the conjugate of a convex function defined in terms of a function of four variables and two affine maps.   Theorem \BDQUADthm\ leads easily to Corollary \BDBIBIcor.   We call Corollary \BDBIBIcor\ a \slant bibivariate existence theorem \endslant\ because it uses two pairs of two spaces. In it, we show how we can compute the conjugate of a convex function defined in terms of two functions of two variables and two affine maps.      We give three consequences of Corollary \BDBIBIcor: Corollaries \BDSZcor--\BDINDSPECcor.   In Corollary \BDSZcor, we show how we can compute the conjugate of a convex function defined in terms of a partial inf--convolution of two functions, and in Corollaries \BDINDcor\ and \BDINDSPECcor, we show how how we can compute the conjugate of a convex function defined in terms of a given convex function and two affine maps.   All through the analysis that we have discussed so far, the conclusion is the existence of a \slant linear\endslant\ functional satisfying certain properties.   Corollary \BDINDSPECcor\ is the first place in which we assume that one of the given maps is \slant linear\endslant\ (rather than \slant affine\endslant).  The statement of Corollary \BDINDSPECcor\ is also interesting in that it does not contain ``$+$'' or ``$-$''.
\par
Section \FNORMsec\ is about (not necessarily locally convex) metrizable linear topological spaces.   Lemma \CAUCHYlem\ and Lemma \FNORMlem\ are technical results, the second of which uses Baire's theorem.   They lead to Theorem \RSthm, which implies that, under the appropriate circumstances, the interiority condition (\TRItwo) is equivalent to a much simpler condition.   The statement of Theorem \RSthm\ is disarmingly simple given how much effort seems to be required to establish it.   The automatic interiority result of Corollary \RScor\ is immediate from Theorem \RSthm.
\par
In Section \COMPLsec, we give applications of Corollary \RScor\ to the existence of linear functionals.   The trivariate existence theorem, Theorem \TRIthm\ is the analog of Theorem \BDTRIthm, and the remaining results in Section \COMPLsec\ follow from Theorem \TRIthm\ in much the same way that the results in Section \BDsec\ followed from Theorem \BDTRIthm.   Corollary \FENcor\ extends the version of the Fenchel duality theorem due to Robinson and Attouch--Brezis that we have already mentioned.   Theorem \QUADthm\ is a second quadrivariate existence theorem, Corollary \BIBIcor\ is a second bibivariate existence theorem that extends a result of Simons, and the result on partial inf--convolutions that appears in Corollary \SZcor\ extends a result of     
Simons--Z\u{a}linescu.   In Corollaries \EEAcor\ and \EEBcor, we apply Corollary \BIBIcor\ to situations where the functions are defined on spaces of the form $E \times E^*$ and $F \times F^*$, where $E$ and $F$ are Banach spaces.   Similar results are true in the context of Section \BDsec, but they are less interesting since, in the situation in which these results are applied, the functions concerned are the Fitzpatrick functions of monotone multifunctions, which are known to be lower semicontinuous. Corollaries \INDcor\ and \INDSPECcor\ also extend results that have been used recently in the study of maximally monotone multifunctions on nonreflexive Banach space.   We refer the reader to \cite\QUADARCHIV,\QUADARCHIVTWO\endcite\ for more details of these applications.
\par
We would like to express our sincere thanks to Constantin Z\u{a}linescu for reading through the first version of this paper, and making a number of suggestions that have improved the exposition enormously. 
\defSection \HBLsec
\medbreak
\centerline{\bf \HBLsec.\quad The existence of linear functionals}
\medskip
\noindent
All vector spaces in this paper are \slant real\endslant.
\defDefinition \PRELIMdef
\medbreak
\noindent
{\bf Definition \PRELIMdef.}\enspace We shall say that $Z$ is a \slant convex combination space (CC space)\endslant\ if $Z$ is a nonempty convex subset of a vector space.   Let $Z$ and $X$ be CC spaces and $B\colon Z \to X$.    We say that $B$ is \slant affine\endslant\ if, for all $x, y \in Z$ and $\lambda \in \,]0,1[$, $B\(\lambda x + (1 - \lambda)y\) = \lambda Bx + (1 - \lambda)By$, and we write $\aff(Z,X)$ for the set off affine functions from $Z$ into $X$.   We write $Z^\flat$ for the set of affine functions from $Z$ into $\r$, and so $Z^\flat = \aff(Z,\r)$.   If $Z$ is a CC space, we write $\PC(Z)$ for the set of all convex functions $k\colon\ Z \to \rbar$ such that $\dom\,k \ne \emptyset$, where $\dom\,k$, the \slant effective domain \endslant of $k$, is defined by
$$\dom\,k := \big\{z \in Z\colon\ k(z) \in \r\big\}.$$
(The ``${\cal P}$'' stands for ``proper'', which is the adjective
frequently used to denote the fact that a function is finite at at
least one point.)   If $X$ is a vector space, we write $X^\prime$ for the set of linear functionals on $X$, the algebraic dual of $X$.
\medskip
The main results of Section \HBLsec\ are the \slant asymmetric sandwich theorem\endslant\ of Theorem \ASANDthm\ and the \slant sublevel set theorem\endslant\ of Theorem \SUBLEVthm. In order to justify this nomenclature for Theorem \ASANDthm, we state K\"onig's original result (see [\KONIG, Theorem 1.7, p. 112]), which can obviously be obtained from Theorem \ASANDthm\ by taking $Z = X$ and $B$ to be the identity map.   (\slant Sublinear\endslant\ is defined in Definition \SUBdef.)
\medbreak
\noindent
{\bf \KST.}\enspace\slant Let $X$ be a vector space, $S\colon\ X \to \r$ be sublinear,  $k \in \PC(X)$ and $S \ge -k$ on $X$.   Then there exists $x^\prime \in X^\prime$ such that $S \ge x^\prime \ge -k$ on $X$.\endslant
\medskip
This theorem is \slant symmetric\endslant\ because $S$ and $k$ are defined on the same set.   By contrast, Theorem \ASANDthm\ is \slant asymmetric\endslant\ because $S$ and $k$ are defined on the different sets $X$ and $Z$.  
\defDefinition \SUBdef
\medbreak
\noindent
{\bf Definition \SUBdef.}\enspace Let $X$ be a nontrivial vector space.   We say that $S\colon\ X \to \r$ is \slant sublinear\endslant\ if
$$S\ \hbox{is \slant subadditive\endslant:}\qquad x_1, x_2 \in X \qlr S(x_1 + x_2) \le S(x_1) + S(x_2)$$
and
$$S\ \hbox{is \slant positively homogeneous\endslant:}\qquad x \in X\ \and\ \lambda > 0 \qlr S(\lambda x) = \lambda S(x).$$
We note that it follows automatically that $S(0) = 0$.
Our results depend on the classical Hahn--Banach theorem for sublinear functionals, which we now state:
\defLemma \HBlem
\medbreak
\noindent
{\bf Lemma \HBlem.}\enspace\slant Let $X$ be a vector
space and $S\colon\ X \to \r$ be sublinear.   Then there exists $x^\prime \in X^\prime$ such that $x^\prime \le S$ on $E$.\endslant
\Proof See Kelly--Namioka, [\KN, 3.4, p.\ 21] for a proof using cones, Rudin, [\RUDIN, Theorem 3.2, p.\ 56--57] for a proof using an extension by subspaces argument, and K\"onig, [\KONIG] and Simons, [\MSF] for a proof using an ordering on sublinear functionals.\qed
\medskip
We now come to the {\bf asymmetric sandwich theorem}.   Remark \CHOICErem\ contains several comments on this result.
\defTheorem \ASANDthm
\medbreak
\noindent
{\bf Theorem \ASANDthm.}\enspace\slant Let $X$ be a vector
space, $S\colon\ X \to \r$ be sublinear, $Z$ be a CC space,  $k \in \PC(Z)$,  $B \in \aff(Z,X)$ and $SB \ge -k$ on $Z$.   Then there exists $x^\prime \in X^\prime$ such that $x^\prime \le S$ on $X$ and $x^\prime B \ge -k$ on $Z$.\endslant
\Proof For all $x \in X$, let
$$T(x) := \inf\nolimits_{z \in Z,\ \lambda > 0}\[S(x + \lambda Bz) + \lambda k(z)\] \in [-\infty,\infty].\meqno\AFFSANDone$$
If $x \in X$, $z \in Z$ and $\lambda > 0$ then
$$S(x + \lambda Bz) + \lambda k(z)\ge S(x + \lambda Bz) - S(\lambda Bz) \ge -S(-x).$$
Taking the infimum over $z \in Z$ and $\lambda > 0$, $T(x) \ge -S(-x) > -\infty$.   On the other hand, fix $z \in \dom\,k$.  Let $x$ be an arbitrary element of $X$.  Then, for all $\lambda > 0$, $T(x) \le S(x + \lambda Bz) + \lambda k(z) \le S(x) + \lambda S(Bz) + \lambda k(z)$.   Letting $\lambda \to 0$, $T(x) \le S(x)$.  Thus
$$T\colon\ X \to \r\quand T \le S\ \on\ X.$$
We now show that $T$ is subadditive.   To this end, let $x_1,\,x_2 \in X$. Let $z_1,\,z_2 \in Z$ and $\lambda_1,\,\lambda_2 > 0$ be arbitrary.   Write $x := x_1 + x_2$, and $z := (\lambda_1z_1 + \lambda_2z_2)/(\lambda_1 + \lambda_2)$.  Then, since $\lambda_1Bz_1 + \lambda_2Bz_2 = (\lambda_1 + \lambda_2)Bz$ and $\lambda_1k(z_1) + \lambda_2k(z_2) \ge (\lambda_1 + \lambda_2)k(z)$,
$$\eqalignno{
\[S(x_1 + \lambda_1Bz_1) + \lambda_1k(z_1)\] &+
\[S(x_2 + \lambda_2Bz_2) + \lambda_2k(z_2)\]\cr
&\ge\textstyle S(x + \lambda_1Bz_1 + \lambda_2Bz_2) +
\lambda_1k(z_1) + \lambda_2k(z_2)\cr
&\ge\textstyle S(x + (\lambda_1 + \lambda_2)Bz) + (\lambda_1 + \lambda_2)k(z) \ge T(x) = T(x_1 + x_2).}$$
Taking the infimum over $z_1$, $z_2$, $\lambda_1$ and $\lambda_2$ gives $T(x_1) + T(x_2) \ge T(x_1 + x_2)$.   Thus $T$ is subadditive.   It is easy to check that $T$ is positively homogeneous, and so $T$ is sublinear.  From Lemma \HBlem, there exists $x^\prime \in X^\prime$ such that $x^\prime \le T$ on $X$.   Since $T \le S$ on $X$, $x^\prime \le S$ on $X$, as required. Finally, let $z$ be an arbitrary element of $Z$.   Taking $\lambda = 1$ in (\AFFSANDone), $k(z) = S({-}Bz + Bz) + k(z) \ge T({-}Bz) \ge x^\prime({-}Bz) = -(x^\prime B)z$, hence $x^\prime B \ge -k$ on X.  This completes the proof of Theorem \ASANDthm.\qed
\defRemark \CHOICErem
\medbreak
\noindent
{\bf Remark \CHOICErem.}\enspace It is worth pointing out that the
definition of the \slant auxiliary sublinear\break functional\endslant, $T$, used to prove Theorem \ASANDthm\ is ``forced'' in the sense that if $x^\prime \in X^\prime$, $x^\prime \le S$ on $X$ and $x^\prime B \ge -k$ on $Z$ then, as the reader can easily verify, $x^\prime \le T$ on $X$.\par
It is easy to see that Theorem \ASANDthm\ follows from the \slant Hahn--Banach--Lagrange theorem\endslant\ of \cite\HBL, Theorem 2.9, p.\ 153\endcite\ or \cite\HBM, Theorem 1.11, p.\ 21\endcite.   On the other hand, Theorem \ASANDthm\ implies the Mazur-Orlicz theorem of \cite\HBL, Lemma 2.4, p.\ 152\endcite\ or \cite\HBM, Lemma 1.6, p.\ 19\endcite, which in turn implies the Hahn--Banach--Lagrange theorem.
\par
While we have presented Theorem \ASANDthm\ as a fairly direct consequence of the Hahn--Banach theorem, one can also establish it using an appropriate version of the Fenchel duality theorem.   We will return to this issue in Remark \SANDFENrem. 
\medskip
We now define the \slant sublevel sets\endslant, and also the \slant conjugate with respect to a real affine function\endslant, of a proper convex function on a CC space.
\defDefinition \STARdef
\medbreak
\noindent
{\bf Definition \STARdef.}\enspace Let $Z$ be a CC space, $\Phi \in \PC(Z)$, $\gamma \in \r$ and $z^\flat \in Z^\flat$.   Then we write $\sigma_\Phi(\gamma)$ for the sublevel set $\big\{z \in Z\colon\ \Phi z < \gamma\big\}$.  The set $\sigma_\Phi(\gamma)$ is convex.   We define $\Phi^*(z^\flat) := \sup_Z[z^\flat - \Phi] \in \rbar$.
\medskip
The next result is the {\bf sublevel set theorem}.   If $X$ is a locally convex space, we write $X^*$ for the set of continuous linear functionals on $X$, the topological dual of $X$.
\defTheorem \SUBLEVthm
\medbreak
\noindent
{\bf Theorem \SUBLEVthm.}\enspace\slant Let $Z$ be a CC space, $X$ be a locally convex space, $B \in \aff(Z,X)$, $\Phi \in \PC(Z)$,
$$Y := \ts\bigcupn_{\lambda > 0}\lambda B(\dom\,\Phi)\ \hbox{be a linear subspace of}\ X,\meqno\ZXone$$
and suppose that there exists $\gamma \in \r$ such that
$$0 \in \intr_YB\(\sigma_\Phi(\gamma)\).\meqno\ZXtwo$$
Then
$$\minn_{x^* \in X^*}\Phi^*(x^*B) = -\inf\Phi\(B^{-1}\{0\}\).\meqno\ZXthree$$\endslant
\Proo From (\ZXtwo), there exists $z_0 \in B^{-1}\{0\} \cap \sigma_\Phi(\gamma)$.  Then $\inf\Phi\(B^{-1}\{0\}\) \le \Phi z_0 < \gamma$, and so
$$\inf\Phi\(B^{-1}\{0\}\) < \gamma < \infty.\meqno\ZXnew$$
Let $x^* \in X^*$ and $z \in B^{-1}\{0\}$.   Then
$$\Phi^*(x^*B)\ge x^*B(z) - \Phi z = x^*(Bz) - \Phi z  = 0 - \Phi z = - \Phi z,$$%
and so\quad $\Phi^*(x^*B) \ge \sup\[-\Phi\(B^{-1}\{0\}\)\] = - \inf\Phi\(B^{-1}\{0\}\)$.\quad   So what we must prove  for (\ZXthree) is that
$$\exs\ x^* \in X^*\ \st\ \Phi^*(x^*B) \le - \inf\Phi\(B^{-1}\{0\}\).\meqno\ZXfour$$
If $\inf\Phi\(B^{-1}\{0\}\) = -\infty$, the result is obvious with $x^* := 0$ so, using (\ZXnew), we can and will suppose that $\inf\Phi\(B^{-1}\{0\}\) \in \r$.   Define $k \in \PC(Z)$ by\quad $k := \Phi - \inf\Phi\(B^{-1}\{0\}\)$.\quad Since $\dom\,k = \dom\,\Phi$, (\ZXone) implies that\quad $Y = \ts\bigcupn_{\lambda > 0}\lambda B(\dom\,k)$.\quad Let $\eta := \gamma - \inf\Phi\(B^{-1}\{0\}\)$. From (\ZXnew) and (\ZXtwo), $\eta > 0$ and there exists a continuous seminorm $S$ on $X$ such that
$$\big\{y \in Y\colon\ Sy < 1\big\} \subset B\(\sigma_k(\eta)\).\meqno\ZXfive$$
From the definition of $k$,
$$z \in B^{-1}\{0\} \qlr k(z) \ge 0.\meqno\ZXsix$$
We now prove that
$$\eta SB \ge -k\ \on\ Z.\meqno\ZXseven$$
To this end, first let $z \in \dom\,k$.  Let  $\mu > S(Bz) \ge 0$.   Then $-Bz/\mu \in Y$ and $S(-Bz/\mu) < 1$, and so (\ZXfive) provides $\zeta \in \sigma_k(\eta)$ such that $-Bz/\mu = B\zeta$,
from which\quad $B\((\mu \zeta  + z)/(\mu + 1)\) = 0$.\quad   Thus, using (\ZXsix) and the convexity of $k$,
$$0 \le k\((\mu \zeta  + z)/(\mu + 1)\) \le \(\mu k(\zeta) + k(z)\)/(\mu + 1) < \(\mu\eta + k(z)\)/(\mu + 1).$$
Letting $\mu \to  S(Bz)$, we see that\quad $0 \le \eta S(Bz) + k(z) = (\eta SB + k)(z)$.\quad   Since this is trivially true if\quad $z \in Z \setminus \dom\,k$,\quad we have established (\ZXseven).  From Theorem \ASANDthm, there exists $x^\prime \in X^\prime$ such that\quad $x^\prime \le \eta S$ on $X$\quad and\quad $x^\prime B \ge -k$ on $Z$.\quad Now any linear functional dominated by $\eta S$ is continuous and so, writing $x^* = -x^\prime$,\quad $x^*B - k \le 0$ on $Z$,\quad that is to say,\quad $x^*B - \Phi \le -\inf\Phi\(B^{-1}\{0\}\)$ on $Z$.\quad   (\ZXfour) follows easily from this.\qed
\defRemark \SUBFUNDrem
\medbreak
\noindent
{\bf Remark \SUBFUNDrem.}\enspace In this remark, we compare Theorem \SUBLEVthm\ with the \slant fundamental duality formula\endslant\ of Z\u{a}linescu, \cite\ZBOOK, Theorem 2.7.1(i), pp. 113--114\endcite.   Let $W$ and $X$ be locally convex spaces, $\Phi \in \PC(W \times X)$, $\pi_X\colon\ W \times X \to X$ be defined by $\pi_X(w,x) := x$ and $\pi_X(\dom\,\Phi) \ni 0$.   Let $Y$ be the linear span of $\pi_X(\dom\,\Phi)$, and suppose that there exists $\gamma \in \r$ such that $0 \in \intr_Y\pi_X\(\sigma_\Phi(\gamma)\)$.   Then it is easily seen that the conditions of Theorem \SUBLEVthm\ are satisfied with $Z := W \times X$ and $B := \pi_X$.   Now, $\pi_X^{-1}\{0\} = W \times \{0\}$ and, for all $x^* \in X^*$, $x^*\pi_X = (0,x^*)\in Z^*$.   Thus Theorem \SUBLEVthm\ implies that $\minn_{x^* \in X^*}\Phi^*(0,x^*) = -\inf\Phi\(W \times \{0\}\)$, which is exactly the conclusion of \cite\ZBOOK, Theorem 2.7.1(i)\endcite.   We now consider the reverse question of deducing Theorem \SUBLEVthm\ from \cite\ZBOOK, Theorem 2.7.1(i)\endcite.   Suppose first that $W$ and $X$ are locally convex spaces, $B \in \aff(W,X)$, $\Phi \in \PC(W)$, and (\ZXone) and (\ZXtwo) are satisfied.   Define $\Psi \in \PC(W \times X)$ by

$$\Psi(w,x) = \cases{\Phi(w)&$(x = Bw);$\cr\infty&$(x \ne Bw)$.}$$
Then $\pi_X(\dom\,\Psi) = B(\dom\,\Phi)$ and, if $\gamma \in \r$, $\pi_X\(\sigma_\Psi(\gamma)\) = B\(\sigma_\Phi(\gamma)\)$.   But then, for all $x^* \in X^*$, $\Psi^*(0,x^*) = \Phi^*(x^*B)$ and $W \times \{0\} = B^{-1}\{0\}$, and so (\ZXthree) follows from \cite\ZBOOK, Theorem 2.7.1(i)\endcite.   This establishes Theorem \SUBLEVthm\ in the special case when $Z$ is a locally convex space.   The general case when $Z$ is a CC space can be deduced from the special case by a series of translations and extensions and using the finest locally convex topology.  
\medskip
Corollary \TRIcor\ is our first {\bf trivariate existence theorem}, in which the function $h$ is defined as the infimum of $\Psi$ over a variable set.   Corollary \TRIcor\ will be applied in Theorems \BDTRIthm\ and \TRIthm.
\defCorollary \TRIcor
\medbreak
\noindent
{\bf Corollary \TRIcor.}\enspace\slant Let $Z$ and $P$ be CC spaces, $X$ be a locally convex space, $B \in \aff(Z,X)$, $A \in \aff(Z,P)$ and $\Psi \in \PC(Z)$.   For all $p \in P$, let
$$h(p) := \ts\inf\Psi\(A^{-1}\{p\} \cap B^{-1}\{0\}\) > -\infty$$
and
$$Y := \ts\bigcupn_{\lambda > 0}\lambda B(\dom\,\Psi)\ \hbox{be a linear subspace of}\ X.\meqno\TRIone$$
Let $p^\flat \in P^\flat$ and $\Phi := \Psi - p^\flat A  \in \PC(Z)$, and suppose that there exists $\gamma \in \r$ such that
$$0 \in \intr_YB\(\sigma_\Phi(\gamma)\).\meqno\TRItwo$$
Then
$$h^*(p^\flat) = \ts\min_{x^* \in X^*}\Psi^*(p^\flat A + x^*B).$$\endslant
\Proof Clearly, $\dom\,\Phi = \dom\,\Psi$, and so (\ZXone) follows from (\TRIone).   Of course, (\ZXtwo) is identical with (\TRItwo).   The result now follows from Theorem \SUBLEVthm\ since   
$$\eqalign{
h^*(p^\flat)
&= \sup\big\{p^\flat(p) - \Psi z\colon\ p \in P,\ z \in A^{-1}\{p\} \cap B^{-1}\{0\}\big\}\cr
&= \sup\big\{p^\flat Az - \Psi z\colon\ p \in P,\ z \in A^{-1}\{p\} \cap B^{-1}\{0\}\big\}\cr
&= \sup\big\{p^\flat Az - \Psi z\colon\ z \in B^{-1}\{0\}\big\} =\sup\[-\Phi\(B^{-1}\{0\}\)\] = -\inf\Phi\(B^{-1}\{0\}\)}$$
and, for all $x^* \in X^*$,
$$\Phi^*(x^*B) = \supn_{z \in Z}\[x^*Bz - \Phi z\] = \supn_{z \in Z}\[p^\flat Az + x^*Bz - \Psi z\] = \Psi^*(p^\flat A + x^*B).\eqno\qed$$\par
\defSection \BDsec
\smallbreak
\centerline{\bf \BDsec.\quad Results with a boundedness hypothesis}
\medskip
\noindent
Theorem \BDTRIthm\ is our second {\bf trivariate existence theorem}, which should be compared with Theorem \TRIthm.   There is an important difference between Corollary \TRIcor\ and Theorem \BDTRIthm.   In Corollary \TRIcor, the choice of the bound $\gamma$ will normally depend on $p^\flat$, while in Theorem \BDTRIthm\ the choice of the bound $\delta$ can be made independently of $p^\flat$.  Theorem \BDTRIthm\ will be used in Corollary \BDFENcor\ and Theorem \BDQUADthm. 
\defTheorem \BDTRIthm
\medbreak
\noindent
{\bf Theorem \BDTRIthm.}\enspace\slant Let $Z$ and $P$ be CC spaces, $X$ be a locally convex space, $B \in \aff(Z,X)$, $A \in \aff(Z,P)$, $\Psi \in \PC(Z)$ and, for all $p \in P$,
$$h(p) := \ts\inf\Psi\(A^{-1}\{p\} \cap B^{-1}\{0\}\) > -\infty.\meqno\BDZPXtwo$$
Suppose that there exist $z_0 \in Z$ and $\delta \in \r$ such that
$$0 \in \intr_XB\(A^{-1}\{Az_0\} \cap \sigma_\Psi(\delta)\).\meqno\BDZPXfour$$
Then 
$$p^\flat \in P^\flat \qlr h^*(p^\flat) = \ts\min_{x^* \in X^*}\Psi^*(p^\flat A + x^*B).\meqno\BDZPXthree$$\endslant
\Proo Let $x \in X$.   If $\lambda$ is sufficiently large then $x/\lambda \in B(\dom\,\Psi)$, from which $x \in \lambda B(\dom\,\Psi)$.   Thus $\ts\bigcupn_{\lambda > 0}\lambda B(\dom\,\Psi) = X$, and (\TRIone) is satisfied.    Furthermore, if\quad $p^\flat \in P^\flat$,\quad $\Phi := \Psi - p^\flat A  \in \PC(Z)$,\quad and\quad $z \in A^{-1}\{Az_0\} \cap \sigma_\Psi(\delta)$\quad then, writing\quad $\gamma:= \delta  - p^\flat Az_0$, 
$$\Phi z = \Psi z - p^\flat Az = \Psi z - p^\flat Az_0 < \delta - p^\flat Az_0 = \gamma.$$
and so\quad $A^{-1}\{Az_0\} \cap \sigma_\Psi(\delta) \subset \sigma_\Phi(\gamma)$.   Consequently, (\BDZPXfour) gives (\TRItwo),   and the result follows from Corollary \TRIcor.\qed
\medbreak
In our first result on {\bf affine Fenchel duality}, Corollary \BDFENcor, which should be compared with Corollary \FENcor, we show how Theorem \BDTRIthm\ leads to a result on the conjugate of a generalized sum of convex functions.   These results can also be deduced from the more general results that follow from Theorem \BDQUADthm\ ---  we have included them here because they provide a model for the somewhat more complex proof of Theorem \BDQUADthm.   Corollary \BDFENcor\ generalizes the classical result of Rockafellar \cite\RTRFENCHEL, Theorem 1\endcite.   In what follows, the product of CC spaces is understood to have the pointwise definition of the convex operation.
\defCorollary \BDFENcor
\medbreak
\noindent
{\bf Corollary \BDFENcor.}\enspace\slant Let $P$ be a CC space, $X$ be a locally convex space, $C \in \aff(P,X)$, $f \in \PC(P)$ and $g \in \PC(X)$.   Suppose that there  exists $p_0 \in \dom\,f$ such that $g$ is finitely bounded above in a neighborhood of $Cp_0$.   Then
$$p^\flat\in P^\flat \qlr (f + gC)^*(p^\flat) = \ts\min_{x^* \in X^*}\[f^*(p^\flat - x^*C) + g^*(x^*)\].\meqno\BDFENtwo$$\endslant
\Proo Let $Z := P \times X$, and define $B \in \aff(Z,X)$ by\quad $B(p,x) :=  x - Cp$,\quad $A \in \aff(Z,P)$ by\quad $A(p,x) := p$,\quad and $\Psi \in \PC(Z)$ by\quad $\Psi(p,x) := f(p) + g(x)$.\quad   If $p \in P$ then it is easy to see that\quad $\Psi\(A^{-1}\{p\} \cap B^{-1}\{0\}\)$\quad is the singleton\quad $\{(f + gC)(p)\}$\quad thus, in the notation of (\BDZPXtwo),\quad $h = f + gC$,\quad hence\quad $h^* = (f + gC)^*$.\quad   Now let $p^\flat\in P^\flat$: then $p^\flat A + x^*B = (p^\flat - x^*C,x^*) \in Z^\flat$.   By direct computation, for all $q^\flat \in P^\flat$, $\Psi^*(q^\flat,x^*) = f^*(q^\flat) + g^*(x^*)$,\quad so the formula for\quad $(f + gC)^*$\quad given in (\BDFENtwo) reduces to the formula for $h^*$ given in (\BDZPXthree). Let $z_0 = (p_0,Cp_0) \in Z$.   By hypothesis, there exists $\gamma \in \r$ such that if $y \in X$ is sufficiently small then\quad $g(Cp_0 + y) < \gamma$,\quad from which\quad $\Psi(p_0,Cp_0 + y) = f(p_0) + g(Cp_0 + y) < f(p_0) + \gamma$.\break   Since\quad $(p_0,Cp_0 + y) \in A^{-1}\{Az_0\}$\quad and\quad $B(p_0,Cp_0 + y) = y$,\quad  (\BDZPXfour) is satisfied with\break $\delta := f(p_0) + \gamma$,\quad and the result follows from Theorem \BDTRIthm.\qed
\defRemark \SANDFENrem
\medbreak
\noindent
{\bf Remark \SANDFENrem.}\enspace Let $X$ be a vector space and $\T$ be the finest locally convex topology on $X$.   Then every sublinear functional on $X$ is finitely bounded above in a neighborhood of every element of $X$.   Thus we can apply Corollary \BDFENcor\ with $P$, $C$, $f$ and $g$ replaced by $Z$, $B$, $k$ and $S$ (respectively), and $p^\flat := 0$, and obtain: \slant Let $S\colon\ X \to \r$ be sublinear, $Z$ be a CC space, $B \in \aff(Z,X)$ and $k \in \PC(Z)$.  Then $(k + SB)^*(0) = \ts\min_{x^* \in (X,\T)^*}\[k^*(-x^*B) + S^*(x^*)\].$\endslant\
Suppose now that $SB \ge -k$ on $Z$.   Then $(k + SB)^*(0) \ge 0$, and so there exists $x^\prime \in X^\prime$ such that $k^*(-x^\prime B) + S^*(x^\prime) \le 0$.   In particular, $S^*(x^\prime) < \infty$, and since $x^\prime - S$ is positively homogeneous, it follows that $x^\prime \le S$ on $X$ and $S^*(x^\prime) = 0$.   Thus $k^*(-x^\prime B) \le 0$, from which $x^\prime B \ge -k$ on $Z$.   Thus, in this somewhat indirect fashion, we obtain Theorem \ASANDthm. 
\medskip
We now come to Theorem \BDQUADthm, our first {\bf quadrivariate existence theorem}, which should be compared with Theorem \QUADthm.   We have taken $V$ to be a vector space because of the expression ``$v - Du$'' that appears in (\BDQUADtwo).   The remaining results in this section all follow easily from Theorem \BDQUADthm.
\defTheorem \BDQUADthm
\medbreak
\noindent
{\bf Theorem \BDQUADthm.}\enspace\slant Let $U$ and $W$ be CC spaces, $V$ be a vector space, $X$ be a locally convex space, $C \in \aff(W,X)$, $D \in \aff(U,V)$, $\Psi \in \PC(U \times V \times W \times X)$ and, for all $(w,v) \in W \times V$,
$$h(w,v) := \ts\inf_{u \in U}\Psi(u,v - Du,w,Cw) > -\infty.\meqno\BDQUADtwo$$
Suppose that there exists $(u_0,v_0,w_0) \in U \times V \times W$ such that $\Psi(u_0,v_0,w_0,\cdot)$ is finitely bounded above in a neighborhood of $Cw_0$.   Then \(with $(w^\flat,v^\flat)(w,v) := w^\flat w + v^\flat v$\), 
$$(w^\flat,v^\flat) \in W^\flat \times V^\flat \qlr h^*(w^\flat,v^\flat) = \ts\min_{x^* \in X^*}\Psi^*(v^\flat D,v^\flat,w^\flat - x^*C,x^*).\meqno\BDQUADthree$$\endslant
\Proo Let $Z :=  U \times V \times W \times X$ and $P := W \times V$, and define $B \in \aff(Z,X)$ by $B(u,v,w,x) := x - Cw$ and $A \in \aff(Z,P)$ by $A(u,v,w,x) := (w,v + Du)$.   If now $(w,v) \in W \times V$ then we have\quad $A^{-1}\{(w,v)\} \cap B^{-1}\{0\} = \big\{(u,v - Du,w,Cw)\colon\ u \in U\big\}$,\break thus the definition of $h$ given in (\BDQUADtwo) reduces to the definition of $h$ given in (\BDZPXtwo).   Now let $(w^\flat,v^\flat) \in W^\flat \times V^\flat$: then\quad $(w^\flat,v^\flat)A + x^*B = (v^\flat D,v^\flat,w^\flat - x^*C,x^*) \in U^\flat \times V^\flat \times W^\flat \times X^*$,\quad and so the formula for $h^*$ given in (\BDQUADthree) reduces to the formula for $h^*$ given in (\BDZPXthree).   Let $z_0 = (u_0,v_0,w_0,Cw_0) \in Z$.   By hypothesis, there exists $\delta \in \r$ such that if $y \in X$ is sufficiently small then\quad $\Psi(u_0,v_0,w_0,Cw_0 + y) < \delta$.\quad   Since $(u_0,v_0,w_0,Cw_0 + y) \in A^{-1}\{Az_0\}$ and\break $B(u_0,v_0,w_0,Cw_0 + y) = y$, (\BDZPXfour) is satisfied, and the result follows from Theorem \BDTRIthm.\qed
\medbreak
Corollary \BDBIBIcor, our first {\bf bibivariate existence theorem}, should be compared with Corollary \BIBIcor.
\defCorollary \BDBIBIcor
\medbreak
\noindent
{\bf Corollary \BDBIBIcor.}\enspace\slant Let $U$ and $W$ be CC spaces, $V$ be a vector space, $X$ be a locally convex space, $C \in \aff(W,X)$, $D \in \aff(U,V)$, $f \in \PC(W \times V)$, $g \in \PC(X \times U)$ and, for all $(w,v) \in W \times V$,
$$h(w,v) := \infn_{u \in U}\[f(w,v - Du) + g(Cw,u)\] > -\infty.\meqno\BDBIBItwo$$
Define $\pi_W\colon\ W \times V \to W$ by $\pi_W(w,v) := w$ and suppose that there exist $w_0 \in \pi_W\dom\,f$ and $u_0 \in U$ such that $g(\cdot,u_0)$ is finitely bounded above in a neighborhood of $Cw_0$.   Then 
$$(w^\flat,v^\flat) \in W^\flat \times V^\flat \qlr h^*(w^\flat,v^\flat) = \minn_{x^* \in X^*}\[f^*(w^\flat -  x^*C,v^\flat) + g^*(x^*,v^\flat D)\].\meqno\BDBIBIthree$$\endslant
\Proof Let $\Psi(u,v,w,x) := f(w,v) + g(x,u)$.   Then the function $h$ as defined in (\BDBIBItwo) reduces to the function $h$ as defined in (\BDQUADtwo).   Since $\Psi^*(u^\flat,v^\flat,w^\flat,x^*) = f^*(w^\flat,v^\flat) + g^*(x^*,u^\flat)$, the formula for $h^*$ given in (\BDBIBIthree) reduces to the formula for $h^*$ given in (\BDQUADthree).   The result is now immediate from Theorem \BDQUADthm.\qed
\medskip
Corollary \BDSZcor, our first result on {\bf partial inf--convolutions}, which should be compared with Corollary \SZcor, follows from Corollary \BDBIBIcor\ by taking $U = V$, $W = X$, and $C$ and $D$ identity maps. 
\defCorollary \BDSZcor
\medbreak
\noindent
{\bf Corollary \BDSZcor.}\enspace\slant Let $V$ be a vector space, $X$ be a locally convex space, $f,g \in \PC(X \times V)$ and, for all $(x,v) \in X \times V$,
$$h(x,v) := \infn_{u \in V}\[f(x,v - u) + g(x,u)\] > -\infty.$$
Define $\pi_X\colon\ X \times V \to X$ by $\pi_X(x,v) := x$ and  suppose that there exist $x_0 \in \pi_X\dom\,f$ and $v_0 \in V$ such that $g(\cdot,v_0)$ is finitely bounded above in a neighborhood of $x_0$.   Then
$$(x^\flat,v^\flat) \in X^\flat \times V^\flat \qlr h^*(x^\flat,v^\flat) = \minn_{x^* \in X^*}\[f^*(x^\flat -  x^*,v^\flat) + g^*(x^*,v^\flat)\].$$\endslant\par
Corollaries \BDINDcor\ and \BDINDSPECcor\ should be compared with Corollaries \INDcor\ and \INDSPECcor.   The function $f$ defined in     (\BDINDfour) is an \slant indicator function\endslant.   Corollaries \BDINDcor, \BDINDSPECcor, \INDcor\ and \INDSPECcor\ are the only results in this paper that use indicator functions.   The statements of these four results are also interesting in that they do not contain ``$+$'' or ``$-$''. 
   
\defCorollary \BDINDcor
\smallbreak
\noindent
{\bf Corollary \BDINDcor.}\enspace\slant Let $U$ and $W$ be CC spaces, $V$ be a vector space, $X$ be a locally convex space, $C \in \aff(W,X)$, $D \in \aff(U,V)$, $g \in \PC(X \times U)$ and,  for all $(w,v) \in W \times V$,
$$h(w,v) := \inf\big\{g(Cw,u)\colon\ u \in U,\ Du = v\big\} > -\infty.\meqno\BDINDone$$
Suppose that
$$y^\flat \in W^\flat\ \and\ \sup y^\flat(W) < \infty \qlr y^\flat = 0,\meqno\BDINDtwo$$
and there exists $(w_0,u_0) \in W \times U$ such that $g(\cdot,u_0)$ is finitely bounded above in a neighborhood of $Cw_0$.   Then
$$\left.\eqalign{(w^\flat,v^\prime) \in W^\flat \times V^\prime\ &\and\ h^*(w^\flat,v^\prime) < \infty \qlr\cr
&h^*(w^\flat,v^\prime) = \min\big\{g^*\(x^*,v^\prime D\)\colon\ x^* \in X^*,\ x^*C = w^\flat\big\}.}\right\}\meqno\BDINDthree$$\endslant
\Proo Define $f \in \PC(W \times V)$ by
$$f(w,v) = \cases{0&$(v = 0);$\cr\infty&$(v \ne 0)$.}\meqno\BDINDfour$$
The definition of $h$ given in (\BDINDone) clearly reduces to that given in (\BDBIBItwo). Let $(w^\flat,v^\prime) \in W^\flat \times V^\prime$ and $h^*(w^\flat,v^\prime) < \infty$.  Then Corollary \BDBIBIcor\ provides us with $x^* \in X^*$ such that\quad $h^*(w^\flat,v^\prime) =\break f^*(w^\flat -  x^*C,v^\prime) + g^*(x^*,v^\prime D)$,\quad and so \quad $f^*(w^\flat -  x^*C,v^\prime) < \infty$.\quad Since $v^\prime0 = 0$, for all $y^\flat \in W^\flat$,\quad $f^*(y^\flat,v^\prime) = \supn_{(w,v) \in W \times V}\[y^\flat w + v^\prime v - f(w,v)\] = \sup y^\flat(W)$.\quad
Thus (\BDINDtwo) implies that\quad $w^\flat = x^*C$\quad and, since\quad $\pi_W\dom\,f = \pi_W\(W \times \{0\}\) = W$\quad the result follows from Corollary \BDBIBIcor.\qed
\defCorollary \BDINDSPECcor
\smallbreak
\noindent
{\bf Corollary \BDINDSPECcor.}\enspace\slant Let $U$ be a CC space, $W$ and $V$ be vector spaces, $X$ be a locally convex space, $C\colon\ W \to X$ be linear, $D \in \aff(U,V)$, $g \in \PC(X \times U)$ and,  for all $(w,v) \in W \times V$,
$$h(w,v) := \inf\big\{g(Cw,u)\colon\ u \in U,\ Du = v\big\} > -\infty.$$
Suppose that there exists $(w_0,u_0) \in W \times U$ such that $g(\cdot,u_0)$ is finitely bounded above in a neighborhood of $Cw_0$.   Then
$$\eqalign{(w^\prime,v^\prime) \in W^\prime \times V^\prime\ &\and\ h^*(w^\prime,v^\prime) < \infty \qlr\cr
&h^*(w^\prime,v^\prime) = \min\big\{g^*\(x^*,v^\prime D\)\colon\ x^* \in X^*,\ x^*C = w^\prime\big\}.}$$\endslant
\Proo This is immediate from Corollary \BDINDcor, since\quad $w^\prime \in W^\prime \lr w^\prime - x^*C \in W^\prime$,\quad and (\BDINDtwo) is true if $W$ is a vector space and $y^\flat \in W^\prime$.\qed
\defSection \FNORMsec
\medbreak
\centerline{\bf \FNORMsec.\quad $(F)$--normed topological vector spaces}
\medskip
\noindent
If $M$ is a vector space, we say that a function $|\cdot|\colon\ M \to [\,0,\infty)$ is an \slant $(F)$--norm\endslant\ if\quad $|x| = 0 \iff x = 0$,\quad for all $x,y \in M$, $|x + y| \le |x| + |y|$,\quad and,\quad for all $\lambda \in [-1,1]$ and $x \in X$, $|\lambda x| \le |x|$.   If $(M,d)$ is a metrizable topological vector space then it is well known that there exists an $(F)$--norm $|\cdot|$ on $M$ such that the topology induced by $d$ on $M$ is identical to that induced on $M$ by the metric $(x,y) \mapsto |x - y|$.   \(See \cite\KOTHE, p.\ 163\endcite.\)   We will say that $(M,|\cdot|)$ is an \slant $(F)$--normed topological vector space\endslant\ if $M$ is a topological vector space, $|\cdot|$ is an $(F)$--norm on $M$, and the topology of $M$ is identical to that induced on $M$ by the metric $(x,y) \mapsto |x - y|$.   We caution the reader that the sequences in $M$ that are $d$--Cauchy are not necessarily the same as those that are $|\cdot|$--Cauchy.   We say that $(Z,|\cdot|)$ is an \slant $(F)$--normed CC space\endslant\ if there exists an $(F)$--normed topological vector space $(\Zh,\wh{|\cdot|})$ such that $Z$ is a nonempty convex subset of $\Zh$, and $|\cdot|$ is the restriction of $\wh{|\cdot|}$ to $Z$.
\medskip
We will need the following result, which seems to be more delicate than the corresponding result in the normed case.   It is not true, even in the simplest case, that the  map $(\lambda,z) \mapsto \lambda z$ is uniformly continuous: for all $\delta > 0$, $(\lambda + \delta)(\lambda + \delta) - \lambda\lambda \to \infty$ as $\lambda \to 
\infty$.
\defLemma \CAUCHYlem
\medbreak
\noindent
{\bf Lemma \CAUCHYlem.}\enspace\slant Let $(M,|\cdot|)$ be an $(F)$--normed topological vector space, $\{\alpha_n\}_{n \ge 1}$ be a \break Cauchy sequence in $\r$, and $\{z_n\}_{n \ge 1}$ be a Cauchy sequence in $(M,|\cdot|)$.   Then $\{\alpha_nz_n\}_{n \ge 1}$ is Cauchy in $(M,|\cdot|)$.\endslant
\Proof Let $\alpha_0 = \lim_{n \to \infty}\alpha_n$.   Let $\eps > 0$.   Since the map $(\lambda,z) \mapsto \lambda z$ is continuous at $(\alpha_0,0)$, there exists $\delta > 0$ such that
$$|\lambda - \alpha_0| < \delta\ \and\ |z| < \delta \qlr |\lambda z| = |\lambda z - \alpha_00| < \eps/3.$$
Since $\alpha_n \to \alpha_0$ in $\r$ and $\{z_n\}_{n \ge 1}$ is Cauchy in $(M,|\cdot|)$, there exists $n_0 \ge 1$ such that
$$n \ge n_0 \qlr |\alpha_n - \alpha_0| < \delta\ \and\ |z_n - z_{n_0}| < \delta.$$
Since the map $\lambda \mapsto \lambda z_{n_0}$ is continuous at $0$, there exists $n_1 \ge n_0$ such that
$$n \ge n_1 \qlr |(\alpha_n - \alpha_{n_1})z_{n_0}| < \eps/3.$$
Now let $n \ge n_1$.    Then $|\alpha_n - \alpha_0| < \delta$, $|z_n - z_{n_0}| < \delta$, $|\alpha_{n_1} - \alpha_0| < \delta$ and $|z_{n_1} - z_{n_0}| < \delta$.   Thus
$$\eqalign{|\alpha_nz_n &- \alpha_{n_1}z_{n_1}|
= |\alpha_n(z_n - z_{n_0}) - \alpha_{n_1}(z_{n_1} - z_{n_0}) + (\alpha_n - \alpha_{n_1})z_{n_0}|\cr
&\le |\alpha_n(z_n - z_{n_0})| + |\alpha_{n_1}(z_{n_1} - z_{n_0})| + |(\alpha_n - \alpha_{n_1})z_{n_0}| < \eps/3 + \eps/3 + \eps/3 = \eps.}$$
This completes the proof of Lemma \CAUCHYlem.\qed 
\defLemma \FNORMlem
\medbreak
\noindent
{\bf Lemma \FNORMlem.}\enspace\slant Let $(Z,|\cdot|)$ be an $(F)$--normed CC space, $X$ be an $(F)$--normed topological vector space, $B \in \aff(Z,X)$, $\Phi \in \PC(Z)$, $Y := \ts\bigcupn_{\lambda > 0}\lambda B(\dom\,\Phi)$ be a complete\break linear subspace of $X$, $\delta \in \r$ and $\delta > \inf \Phi\(B^{-1}\{0\}\)$.   Then:
\par
\noindent
{\rm(a)}\enspace $Y = \ts\bigcup_{i \ge 1}iB\(\sigma_\Phi(\delta)\)$.
\par
\noindent
{\rm(b)}\enspace If $q, k \ge 1$, let\quad $R(q,k) := \big\{z \in \sigma_\Phi(\delta)\colon\ |(1/k){z}| < 2^{-q}\big\}$.\quad  Let $q \ge 1$.  Then there exists $k \ge 1$ such that $0 \in \intr_Y\overline{B\(R(q,k)\)}$.
\par
\noindent
{\rm(c)}\enspace $0 \in \intr_Y\overline{B\(\sigma_\Phi(\delta)\)}$.\endslant 
\Proof(a)\enspace We can fix $z_0 \in \sigma_\Phi(\delta)$ so that $Bz_0 = 0$.   If $y \in Y$, there exist $\lambda > 0$ and $\zeta \in \dom\,\Phi$ such that $y = \lambda B\zeta$.   If $i \ge 1$ and $i > \lambda$ then, since\quad $\Phi z_0 < \delta$,\quad $\Phi\zeta \in \r$\quad and\quad $\Phi\((1 - \lambda/i)z_0 + (\lambda/i)\zeta\) \le (1 - \lambda/i)\Phi z_0 + (\lambda/i)\Phi\zeta$,\quad we can choose $i$ so large that\quad $\Phi\((1 - \lambda/i)z_0 + (\lambda/i)\zeta\) < \delta$,\quad and so\quad $(1 - \lambda/i)z_0 + (\lambda/i)\zeta \in \sigma_\Phi(\delta)$.\quad   Then\quad 
$$y = \lambda B\zeta = iB\((1 - \lambda/i)z_0 + (\lambda/i)\zeta\) \in iB\(\sigma_\Phi(\delta)\).$$
Since $B\(\sigma_\Phi(\delta)\) \subset B(\dom\,\Phi) \subset Y$, this completes the proof of (a).
\par
(b)\enspace For all $z \in \sigma_\Phi(\delta)$, $|(1/k)z| \to 0$ as $k \to \infty$.   Thus it follows from (a) that $Y = \ts\bigcup_{i,m \ge 1}iB\(R(q + 1,m)\)$, and so Baire's theorem provides us with $i,m \ge 1$ and $y_0 \in Y$ such that $y_0 \in \intr_Yi\overline{B\(R(q + 1,m)\)}$.   Since $-y_0 \in Y$, there exist $j,n \ge 1$ and $z_2 \in R(q + 1,n)$ such that $-y_0 = jBz_2$.   Let $k := m \vee n$, $z_1 \in R(q + 1,m)$, and write $z_3 := (iz_1 + jz_2)/(i + j)$.   Then\quad $im/\((i + j)k\) \le 1$\quad and\quad $jn/\((i + j)k\) \le 1$,\quad and so
$$\eqalign{\ts |z_3/k| &= \ts\big|\[im/\((i + j)k\)\](z_1/m) + \[jn/\((i + j)k\)(z_2/n)\]\big|\cr
&\le \ts\big|\[im/\((i + j)k\)\](z_1/m)\big|  + \big|\[jn/\((i + j)k\)(z_2/n)\]\big|\cr
&\le |z_1/m|  + |z_2/n| < 2^{-q - 1} + 2^{-q - 1} = 2^{-q}.}$$
Consequently, $z_3 \in \(R(q,k)\)$.   Since\quad  $iBz_1 + jBz_2 = (i + j)Bz_3$,\quad we have proved that \quad $iB\(R(q + 1,m)\) + jBz_2 \subset (i + j)B\(R(q,k)\)$.\quad Thus
$$0 = y_0 - y_0 \in \intr_Y\overline{iB\(R(q + 1,m)\)} + jBz_2,$$
from which
$$0 \in \intr_Y\overline{iB\(R(q + 1,m)\) + jBz_2} \subset \intr_Y\overline{(i + j)B\(R(q,k)\)}.$$
This gives (b).
\par
(c)\enspace This is immediate from (b), since $R(q,k) \subset \sigma_\Phi(\delta)$.\qed
\medskip
In the sequel, we write $\caff(Z,X)$ for $\{B \in \aff(Z,X)\colon\ B\ \hbox{is continuous}\}$ and $\PCLSC(Z)$ for $\big\{f \in \PC(Z)\colon\ f\ \hbox{is lower semicontinuous}\big\}$.            
Theorem \RSthm\ below is a considerable sharpening of the result proved in Rodrigues--Simons, \cite\RS, Lemma 1, pp.\ 1072--1073\endcite.   We note that we can make the constant $\alpha$ in (\RSfour) as close as we like to $1$ by increasing the rate of growth of $\{k_j\}_{j \ge 1}$. 
\defTheorem \RSthm
\medbreak
\noindent
{\bf Theorem \RSthm.}\enspace\slant Let $(Z,|\cdot|)$ be a complete $(F)$--normed CC space, $X$ be an $(F)$--normed topological vector space,\quad $B \in \caff(Z,X)$,\quad $\Phi \in \PCLSC(Z)$,\quad $Y := \ts\bigcupn_{\lambda > 0}\lambda B(\dom\,\Phi)$ be a complete linear subspace of $X$, and $\gamma \in \r$.   Then the conditions {\rm(\RSone)--(\RSthree)} are equivalent:
$$0 \in \intr_YB\(\sigma_\Phi(\gamma)\)\meqno\RSone$$
$$0 \in B\(\sigma_\Phi(\gamma)\).\meqno\RStwo$$
$$\gamma > \inf \Phi\(B^{-1}\{0\}\).\meqno\RSthree$$\endslant
\Proof It is immediate that (\RSone)$\lr$(\RStwo).\enspace If (\RStwo) is true then  there exists $z \in \sigma_\Phi(\gamma)$ such that $Bz = 0$, and so (\RSthree) is true.
\par
Suppose, finally, that (\RSthree) is true.   Choose $\delta \in \r$ so that $\inf \Phi\(B^{-1}\{0\}\) < \delta < \gamma$.   Let $k_1 = 1$.   From Lemma \FNORMlem(b), for all $q \ge 2$, there exists $k_q \ge 2^{q - 1}$ such that $0 \in \intr_Y\overline{B\(R(q,k_q)\)}$.   For all $q \ge 1$, let $\eta_q:= 1/k_q$.   Write $\alpha := 1/\sum_{q = 1}^\infty\eta_q \in \,]0,1[\,$.   We will prove that
$$\alpha\overline{B\(\sigma_\Phi(\delta)\)} \subset B\(\sigma_\Phi(\gamma)\).\meqno\RSfour$$
To this end, let $y \in \overline{B\(\sigma_\Phi(\delta)\)} = \overline{\eta_1B\(\sigma_\Phi(\delta)\)}$.   Let $\big\{V_q\}_{q \ge 1}$ be a base for the neighborhoods of $0$ in $X$.      Then there exists $z_1 \in \sigma_\Phi(\delta)$ such that $y - \eta_1Bz_1 \in \overline{\eta_2B\(R(2,k_2)\)} \cap V_1$.   Similarly, there exists $z_2 \in R(2,k_2)$ such that $y - \eta_1Bz_1 - \eta_2Bz_2 \in \overline{\eta_3B\(R(3,k_3)\)} \cap V_2$.
Continuing this process inductively, we end up with a sequence $\{z_q\}_{q \ge 1}$ such that,
$$\all\ q \ge 2,\quad z_q \in R(q,k_q) \subset \sigma_\Phi(\delta)\quad \and,\quad \all\ m \ge 1,\quad y - \ts\sum_{i = 1}^m\eta_iBz_i \in V_m.\meqno\RSfive$$  
In what follows, let $(\Zh,\wh{|\cdot|})$ be an $(F)$--normed topological vector space such that $Z \subset \Zh$ and $|\cdot|$ is the restriction of $\wh{|\cdot|}$ to $Z$.   For all $n \ge 1$, let $s_n := \sum_{q = 1}^n\eta_qz_q \in \Zh$ and $\alpha_n := 1/\sum_{q = 1}^n\eta_q \in \,]0,1[\,$.   For all $q \ge 1$, $z_q \in Z$ and $\Phi(z_q) < \delta$ and so, for all $n \ge 1$, \
$$\alpha_ns_n \in Z \quand \Phi(\alpha_ns_n) < \delta.\meqno\RSsix$$
Let $1 \le m < n$.   Then \(since $z_q \in R(q,k_q)$\) $|s_n - s_m| = |\sum_{q = m + 1}^n\eta_qz_q| \le \sum_{q = m + 1}^n|\eta_qz_q| = \sum_{q = m + 1}^n|z_q/k_q| \le \sum_{q = m + 1}^n2^{-q} < 2^{-m}$.   Thus the sequence $\{s_n\}_{n \ge 1}$ is Cauchy in $\Zh$.   Since the sequence $\{\alpha_n\}_{n \ge 1}$ is convergent in $\r$, Lemma \CAUCHYlem\ implies that the sequence $\{\alpha_ns_n\}_{n \ge 1}$ is Cauchy in $Z$, and so the completeness of $Z$ implies that $\lim_{n \to \infty}\alpha_ns_n$ exists in $Z$.  Write $z_0$ for this limit.   Passing to the limit in (\RSsix) and using the lower semicontinuity of $\Phi$, $\Phi z_0 \le \delta$, from which $z_0 \in \sigma_\Phi(\gamma)$.   Since $B$ is continuous on $Z$, it follows from (\RSfive) and the observation that $\sum_{q = 1}^n\alpha_n\eta_q = 1$ that
$$\eqalign{\alpha y
&= \ts\(\lim_{n \to \infty}\alpha_n\)\(\lim_{n \to \infty}\sum_{q = 1}^n\eta_qBz_q\) = \ts\lim_{n \to \infty}\(\alpha_n\sum_{q = 1}^n\eta_qBz_q\)\cr
&= \ts \lim_{n \to \infty}\(\sum_{q = 1}^n\alpha_n\eta_qBz_q\)
= \ts \lim_{n \to \infty}B\(\sum_{q = 1}^n\alpha_n\eta_qz_q\)
= \ts \lim_{n \to \infty}B(\alpha_ns_n) = Bz_0.}$$
Consequently, $\alpha y \in B\(\sigma_\Phi(\gamma)\)$, which gives (\RSfour).  From Lemma \FNORMlem(c), $0 \in  \intr_Y\overline{B\(\sigma_\Phi(\delta)\)}$, and (\RSone) follows from (\RSfour).\qed
\medskip
The final result in this section is about {\bf automatic interiority}.
\defCorollary \RScor
\medbreak
\noindent
{\bf Corollary \RScor.}\enspace\slant Let $(Z,|\cdot|)$ be a complete $(F)$--normed CC space, $X$ be an $(F)$--normed topological vector space, $B \in \caff(Z,X)$, $\Phi \in \PCLSC(Z)$   and $Y := \ts\bigcupn_{\lambda > 0}\lambda B(\dom\,\Phi)$ be a complete linear subspace of $X$.   Then there exists $\gamma \in \r$ such that $0 \in \intr_YB\(\sigma_\Phi(\gamma)\)$, that is to say, {\rm(\TRItwo)} is satisfied.\endslant
\Proof Fix $z_0 \in \dom\,\Phi$ such that $Bz_0 = 0$, and let $\gamma > \Phi z_0$.   Then $0 \in B\(\sigma_\Phi(\gamma)\)$, and the result follows from Theorem \RSthm\((\RStwo)$\lr$(\RSone)\).\qed
\defSection \COMPLsec
\bigbreak
\centerline{\bf \COMPLsec.\quad Results that use completeness}
\medskip
\noindent
Theorem \TRIthm\ is our third {\bf trivariate existence theorem}, which should be compared with Theorem \BDTRIthm.   We write $Z^\sharp := \{z^\flat \in Z^\flat\colon\ z^\flat\ \hbox{is continuous}\} = \caff(Z,\r)$.   We note that a \slant Fr\'echet space\endslant\ is a complete $(F)$--normed \slant locally convex\endslant\ topological vector space.   In this connection, Banach's definition of space of \slant type (F)\endslant\ \(see \cite\BANACH, p. 35\endcite\) does not require either local convexity, or the continuity of the map $(\lambda,x) \mapsto \lambda x$. 
\defTheorem \TRIthm
\medbreak
\noindent
{\bf Theorem \TRIthm.}\enspace\slant Let $Z$ be a complete $(F)$--normed CC space, $P$ be a CC space, $X$ be a Fr\'echet space, $B \in \caff(Z,X)$, $A \in \aff(Z,P)$, $\Psi \in \PCLSC(Z)$, $\ts\bigcupn_{\lambda > 0}\lambda B(\dom\,\Psi)$ be a closed linear subspace of $X$ and, for all $p \in P$,
$$h(p) := \ts\inf\Psi\(A^{-1}\{p\} \cap B^{-1}\{0\}\) > -\infty.\meqno\ZPXtwo$$
Then
$$p^\flat \in P^\flat\ \and\ p^\flat A \in Z^\sharp \qlr h^*(p^\flat) = \ts\min_{x^* \in X^*}\Psi^*(p^\flat A + x^*B).\meqno\ZPXthree$$\endslant
\Proo Let $p^\flat \in P^\flat$ and $p^\flat A \in Z^\sharp$.   Let $\Phi:= \Psi - p^\flat A \in \PCLSC(Z)$.   Since $\dom\,\Phi = \dom\,\Psi$, $\ts\bigcupn_{\lambda > 0}\lambda B(\dom\,\Phi)$ is a closed subspace of a  Fr\'echet space, and thus complete.   The result now follows from Corollaries \RScor\ and \TRIcor.\qed
\medskip
In Corollary \FENcor\ below, which should be compared with Corollary \BDFENcor, we show how Theorem \TRIthm\ leads to a new result on the conjugate of a generalized sum of convex functions.   Corollary \FENcor\ is a generalization of the generalization of the Robinson--Attouch--Brezis theorem to Fr\'echet spaces that first appeared in Rodrigues--Simons, \cite\RS, Theorem 6, p.\ 1076\endcite.   A result similar to the latter, with slightly more restrictive hypotheses, had been established previously by Z\u{a}linescu in \cite\ZZOR, Corollary 4, p. A91\endcite.   We also refer the reader to Z\u{a}linescu, \cite\ZIA, Corollary 2.2, p.\ 22 and Theorem 4.3, p.\ 32\endcite\ for earlier results in this direction.
\defCorollary \FENcor
\medbreak
\noindent
{\bf Corollary \FENcor.}\enspace\slant  Let $P$ be a complete $(F)$--normed CC space, $X$ be a Fr\'echet space, $C \in \caff(P,X)$, $f \in \PCLSC(P)$, $g \in \PCLSC(X)$, and $\ts\bigcupn_{\lambda > 0}\lambda\[\dom\,g - C(\dom\,f)\]$ be a closed linear subspace of $X$.   Then
$$p^\sharp \in P^\sharp \qlr (f + gC)^*(p^\sharp) = \ts\min_{x^* \in X^*}\[f^*(p^\sharp - x^*C) + g^*(x^*)\].\meqno\FENtwo$$\endslant
\Proo Let $Z := P \times X$, and define $B \in \caff(Z,X)$ by $B(p,x) := x - Cp$, $A \in \caff(Z,P)$ by $A(p,x) := p$, and $\Psi \in \PCLSC(Z)$ by $\Psi(p,x) := f(p) + g(x)$.   In the notation of (\ZPXtwo), $h = f + gC$, hence $h^* = (f + gC)^*$.   Now let $p^\sharp \in P^\sharp$: then $p^\sharp A + x^*B = (p^\sharp - x^*C,x^*) \in Z^\sharp$.  By direct computation, for all $q^\sharp \in P^\sharp$, $\Psi^*(q^\sharp,x^*) = f^*(q^\sharp) + g^*(x^*)$, so the formula for $(f + gC)^*$ given in (\FENtwo) reduces to the formula for $h^*$ given in (\ZPXthree).  Since $\dom\,\Psi = \dom\,f \times \dom\,g$ and $B(\dom\,\Psi) = \dom\,g - C(\dom\,f)$, the result is immediate from Theorem \TRIthm.\qed
\medskip
Theorem \QUADthm\ is our second {\bf quadrivariate existence theorem}, which should be compared with Theorem \BDQUADthm.      The remaining results in this paper all follow from Theorem \QUADthm.
\defTheorem \QUADthm
\medbreak
\noindent
{\bf Theorem \QUADthm.}\enspace\slant Let $U$ and $W$  be complete $(F)$--normed CC spaces, $V$   be a complete $(F)$--normed topological vector space, $X$ be a Fr\'echet space, $C \in \caff(W,X)$, $D \in \caff(U,V)$, $\Psi \in \PCLSC(U \times V \times W \times X)$, $\ts\bigcupn_{\lambda > 0}\lambda\big\{x - Cw\colon\ (u,v,w,x) \in \dom\,\Psi\big\}$ be a closed linear subspace of $X$ and, for all $(w,v) \in W \times V$,
$$h(w,v) := \ts\inf_{u \in U}\Psi(u,v - Du,w,Cw) > -\infty.\meqno\QUADtwo$$
Then 
$$(w^\sharp,v^\sharp) \in W^\sharp \times V^\sharp \qlr h^*(w^\sharp,v^\sharp) = \ts\min_{x^* \in X^*}\Psi^*(v^\sharp D,v^\sharp,w^\sharp - x^*C,x^*).\meqno\QUADthree$$\endslant
\Proo Let $Z :=  U \times V \times W \times X$ and $P := W \times V$, and define $B \in \caff(Z,X)$ by $B(u,v,w,x) := x - Cw$ and $A \in \caff(Z,P)$ by $A(u,v,w,x) := (w,v + Du)$.   If now $(w,v) \in W \times V$ then we have\quad $A^{-1}\{(w,v)\} \cap B^{-1}\{0\} = \big\{(u,v - Du,w,Cw)\colon\ u \in U\big\}$,\break thus the definition of $h$ given in (\QUADtwo) reduces to the definition of $h$ given in (\ZPXtwo).   Now let $(w^\sharp,v^\sharp) \in W^\sharp \times V^\sharp$: then\quad $(w^\sharp,v^\sharp)A + x^*B = (v^\sharp D,v^\sharp,w^\sharp - x^*C,x^*) \in U^\sharp \times V^\sharp \times W^\sharp \times X^*$,\quad so the formula for $h^*$ given in (\QUADthree) reduces to the formula for $h^*$ given in (\ZPXthree).    Since $B(\dom\,\Psi) = \big\{x - Cw\colon\ (u,v,w,x) \in \dom\,\Psi\big\}$, the result follows from Theorem \TRIthm.\qed
\medbreak
Corollary \BIBIcor\ is our second {\bf bibivariate existence theorem}, which should be compared with Corollary \BDBIBIcor. This generalizes the result for Banach spaces that first appeared in Simons, \cite\QUADARCHIV, Theorem 3, pp.\ 2--4\endcite.
\defCorollary \BIBIcor
\medbreak
\noindent
{\bf Corollary \BIBIcor.}\enspace\slant Let $U$ and $W$  be complete $(F)$--normed CC spaces, $V$  be a complete $(F)$--normed topological vector space, $X$ be a Fr\'echet space, $C \in \caff(W,X)$, $D \in \caff(U,V)$, $f \in \PCLSC(W \times V)$, $g \in \PCLSC(X \times U)$ and, for all $(w,v) \in W \times V$,
$$h(w,v) := \infn_{u \in U}\[f(w,v - Du) + g(Cw,u)\] > -\infty.\meqno\BIBIone$$
Define $\pi_X\colon\ X \times U \to X$ and $\pi_W\colon\ W \times V \to W$ by  $\pi_X(x,u) := x$ and $\pi_W(w,v) := w$.   If $\ts\bigcupn_{\lambda > 0}\lambda\[\pi_X\,\dom\,g - C\pi_W\,\dom\,f\]$ is a closed linear subspace of $X$ then
$$(w^\sharp,v^\sharp) \in W^\sharp \times V^\sharp \qlr h^*(w^\sharp,v^\sharp) = \minn_{x^* \in X^*}\[f^*(w^\sharp -  x^*C,v^\sharp) + g^*(x^*,v^\sharp D)\].\meqno\SZthree$$\endslant
\Proo Let $\Psi(u,v,w,x) := f(w,v) + g(x,u)$.   Then the definition of $h$ given in (\BIBIone) reduces to the definition of $h$ given in (\QUADtwo).   Now let $(w^\sharp,v^\sharp) \in W^\sharp \times V^\sharp$: then\quad $\Psi^*(u^\sharp,v^\sharp,w^\sharp,x^*) = f^*(w^\sharp,v^\sharp) + g^*(x^*,u^\sharp)$, so the formula for $h^*$ given in (\SZthree) reduces to the formula for $h^*$ given in (\QUADthree).   Since $\{x - Cw\colon\ (u,v,w,x) \in \dom\,\Psi\big\} = \pi_X\,\dom\,g - C\pi_W\,\dom\,f$, the result now follows from Theorem \QUADthm.\qed
\medskip
Corollary \SZcor, our second result on {\bf partial inf--convolutions}, which should be compared with Corollary \BDSZcor, follows from Corollary \BIBIcor\ by taking $U = V$, $W = X$, and $C$ and $D$ to be identity maps.   Corollary \SZcor\ generalizes the result for Banach spaces that first appeared in Simons--Z\u{a}linescu, \cite\SZNZ, Theorem 4.2, pp.\ 9--10\endcite.  
\defCorollary \SZcor
\medbreak
\noindent
{\bf Corollary \SZcor.}\enspace\slant Let $V$  be a complete $(F)$--normed topological vector space, $X$ be a Fr\'echet space, $f,g \in \PCLSC(X \times V)$ and, for all $(x,v) \in X \times V$,
$$h(x,v) := \infn_{u \in U}\[f(x,v - u) + g(x,u)\] > -\infty.$$
Define $\pi_X\colon\ X \times V \to X$ by $\pi_X(x,v) := x$.   If $\ts\bigcupn_{\lambda > 0}\lambda\[\pi_X\,\dom\,g - \pi_X\,\dom\,f\]$ is a closed linear subspace of $X$ then
$$(x^\sharp,v^\sharp) \in X^\sharp \times V^\sharp \qlr h^*(x^\sharp,v^\sharp) = \minn_{x^* \in X^*}\[f^*(x^\sharp -  x^*,v^\sharp) + g^*(x^*,v^\sharp)\].$$\endslant
\par
Corollary \EEAcor\ is immediate from Corollary \BIBIcor\ with $U = F^*$, $V = E^*$, $W = E$ and $X = F$.   Corollary \EEAcor\ is a slight generalization of the result for Banach spaces that first appeared in \cite\QUADARCHIV, Theorem 5(a),\ pp.\ 4--5\endcite.        
\defCorollary \EEAcor
\medbreak
\noindent
{\bf Corollary \EEAcor.}\enspace\slant Let $E$ and $F$ be Banach spaces, $C \in \caff(E,F)$, $D \in \caff(F^*,E^*)$, $f \in \PCLSC(E \times E^*)$, $g \in \PCLSC(F \times F^*)$ and, for all $(x,x^*) \in E \times E^*$,
$$h(x,x^*) := \infn_{y^* \in F^*}\[f(x,x^* - Dy^*) + g(Cw,y^*)\] > -\infty.$$
Define $\pi_F\colon\ F \times F^* \to F$ and $\pi_E\colon\ E \times E^* \to E$ by  $\pi_F(y,y^*) := y$ and $\pi_E(x,x^*) := x$.   If $\ts\bigcupn_{\lambda > 0}\lambda\[\pi_F\,\dom\,g - C\pi_E\,\dom\,f\]$ is a closed linear subspace of $F$ then
$$(x^\sharp,v^\sharp) \in E^\sharp \times E^{*\sharp} \qlr h^*(x^\sharp,v^\sharp) = \minn_{y^* \in F^*}\[f^*(x^\sharp -  y^*C,v^\sharp) + g^*(y^*,v^\sharp D)\].$$\endslant
\par
Corollary \EEBcor\ is immediate from Corollary \BIBIcor\ with $U = F$, $V = E$, $W = E^*$ and $X = F^*$, and changing the order of the arguments of $f$, $g$ and $h$.   Corollary \EEBcor\ generalizes the result for Banach spaces that first appeared in \cite\QUADARCHIV, Theorem 5(b),\ pp.\ 4--5\endcite.
\defCorollary \EEBcor
\medbreak
\noindent
{\bf Corollary \EEBcor.}\enspace\slant Let $F$ and $E$  be Banach spaces, $C \in \caff(E^*,F^*)$, $D \in \caff(F,E)$, $f \in \PCLSC(E \times E^*)$, $g \in \PCLSC(F \times F^*)$ and, for all $(x,x^*) \in E \times E^*$,
$$h(x,x^*) := \infn_{y \in F}\[f(x - Dy,x^*) + g(y,Cx^*)\] > -\infty.$$
Define $\pi_{F^*}\colon\ F \times F^* \to F^*$ and $\pi_{E^*}\colon\ E \times E^* \to E^*$ by  $\pi_{F^*}(y,y^*) := y^*$ and $\pi_{E^*}(x,x^*) := x^*$.   If $\ts\bigcupn_{\lambda > 0}\lambda\[\pi_{F^*}\,\dom\,g - C\pi_{E^*}\,\dom\,f\]$ is a closed linear subspace of $F^*$ then
$$(x^\sharp,w^\sharp) \in E^\sharp \times {E^*}^\sharp \qlr h^*(x^\sharp,w^\sharp) = \minn_{y\dst \in F\dst}\[f^*(x^\sharp,w^\sharp -  y\dst C) + g^*(x^\sharp D,y\dst)\].$$\endslant
\par
Corollaries \INDcor\ and \INDSPECcor, which should be compared with  Corollaries \BDINDcor\ and \BDINDSPECcor, generalize the result for Banach spaces that first appeared in \cite\QUADARCHIV, Theorem 21, pp.\ 12--13\endcite.   This latter result also appeared in \cite\QUADARCHIVTWO, Theorem 6\endcite, and was taken there as the starting point for proofs of a number of results that have already appeared in this paper, as well as many others which have applications to the theory of strongly representable multifunctions.   We refer the reader to \cite\QUADARCHIVTWO\endcite\ for more details of these applications.
\defCorollary \INDcor
\medbreak
\noindent
{\bf Corollary \INDcor.}\enspace\slant Let $U$ and $W$ be complete  $(F)$--normed CC spaces, $V$ be a complete $(F)$--normed topological vector space, $X$ be a Fr\'echet space, $C \in \caff(W,X)$, $D \in \caff(U,V)$, $g \in \PCLSC(X \times U)$ and, for all $(w,v) \in W \times V$,
$$h(w,v) := \inf\big\{g(Cw,u)\colon\ u \in U,\ Du = v\big\} > -\infty.$$
Suppose that
$$y^\sharp \in W^\sharp\ \and\ \sup y^\sharp(W) < \infty \qlr y^\sharp = 0,$$
and define $\pi_X\colon\ X \times U \to X$ by  $\pi_X(x,u) := x$.   If $\ts\bigcupn_{\lambda > 0}\lambda\[\pi_X\,\dom\,g - C(W)\]$ is a closed linear subspace of $X$ then
$$\eqalign{(w^\sharp,v^*) \in W^\sharp \times V^*\ &\and\ h^*(w^\sharp,v^*) < \infty \qlr\cr
&h^*(w^\sharp,v^*) = \min\big\{g^*\(x^*,v^* D\)\colon\ x^* \in X^*,\ x^*C = w^\sharp\big\}.}$$\endslant
\Proo We define $f \in \PCLSC(W \times V)$ as in (\BDINDfour), and   $\pi_W\colon\ W \times V \to W$ by $\pi_W(w,v) := w$.      Then $\pi_W\,\dom\,f = \pi_W\(W \times \{0\}\) = W$.   The rest of the proof now proceeds exactly as in Corollary \BDINDcor, only using Corollary \BIBIcor\ instead of Corollary \BDBIBIcor.\qed 
\defCorollary \INDSPECcor
\medbreak
\noindent
{\bf Corollary \INDSPECcor.}\enspace\slant Let $U$ be a complete $(F)$--normed CC space, $W$ and $V$ be complete $(F)$--normed topological vector spaces, $X$ be a Fr\'echet space, $C\colon\ W \to X$ be continuous and linear, $D \in \caff(U,V)$, $g \in \PCLSC(X \times U)$ and, for all $(w,v) \in W \times V$,
$$h(w,v) := \inf\big\{g(Cw,u)\colon\ u \in U,\ Du = v\big\} > -\infty,$$
and define $\pi_X\colon\ X \times U \to X$ by  $\pi_X(x,u) := x$.   If $\,\ts\bigcupn_{\lambda > 0}\lambda\[\pi_X\,\dom\,g - C(W)\]$ is a closed linear subspace of $X$ then
$$\eqalign{(w^*,v^*) \in W^* \times V^*\ &\and\ h^*(w^*,v^*) < \infty \qlr\cr
&h^*(w^*,v^*) = \min\big\{g^*\(x^*,v^* D\)\colon\ x^* \in X^*,\ x^*C = w^*\big\}.}$$\endslant
\Proo This follows from Corollary \INDSPECcor\ in exactly the same way that Corollary \BDINDSPECcor\ followed from Corollary \BDINDcor.\qed
\bigbreak
\centerline{\bf References}
\nmbr\BANACH
\item{[\BANACH]}S. Banach, \slant Th\'eorie des op\'erations linŽaires\endslant, Chelsea Publishing Co, New York, (1955).
\nmbr\KN
\item{[\KN]}J. L. Kelley, I. Namioka, and co-authors, \slant Linear 
Topological Spaces\endslant, D. Van Nostrand Co., Inc.,
Princeton -- Toronto -- London -- Melbourne (1963).
\nmbr\KONIG
\item{[\KONIG]}H. K\"onig, \slant Some Basic Theorems in Convex
Analysis\endslant, in ``Optimization and operations research'', edited
by B. Korte, North-Holland (1982).
\nmbr\KOTHE
\item{[\KOTHE]}G. K\"othe, \slant Topological vector spaces I\endslant\ (Translated from the German by D. J. H. Garling), Grundlehren der mathematischen Wissenschaften, Vol 159, Springer-Verlag New York 1969.
\nmbr\RTRFENCHEL
\item{[\RTRFENCHEL]} R. T. Rockafellar, \slant Extension of Fenchel's
duality theorem for convex functions\endslant, Duke Math. J. {\bf33}
(1966),  81--89.
\nmbr\RS
\item{[\RS]}B. Rodrigues and S. Simons, \slant Conjugate functions and subdifferentials in nonnormed situations for operators with complete graphs\endslant, Nonlinear Anal. {\bf 12} (1988), 1069Ð-1078.
\nmbr\RUDIN
\item{[\RUDIN]}W. Rudin, \slant Functional analysis\endslant,
McGraw-Hill, New York (1973).
\nmbr\MSF
\item{[\MSF]}S. Simons, \slant Minimal sublinear functionals\endslant, 
Studia Math. {\bf 37} (1970), 37--56.
\nmbr\HBL
\item{[\HBL]}-----, \slant The Hahn--Banach--Lagrange theorem\endslant, Optimization, {\bf 56} (2007), 149--169.
\nmbr\HBM
\item{[\HBM]}-----, \slant From Hahn--Banach to monotonicity\endslant, 
Lecture Notes in Mathematics, {\bf 1693},\break second edition, (2008), Springer--Verlag.
\nmbr\QUADARCHIV
\item{[\QUADARCHIV]}-----, \slant Quadrivariate versions of the AttouchÐ-Brezis theorem and strong representability\endslant, arXiv:0809.0325, posted September 1, 2008.  
\nmbr\QUADARCHIVTWO
\item{[\QUADARCHIVTWO]}-----, \slant Quadrivariate existence theorems and strong representability\endslant, arXiv:0809.0325v2, posted February 22, 2011.  
\nmbr\SZNZ
\item{[\SZNZ]}S. Simons and C. Z\u{a}linescu, \slant Fenchel duality,
Fitzpatrick functions and maximal monotonicity\endslant, J. of
Nonlinear and  Convex Anal., {\bf 6} (2005), 1--22.
\nmbr\ZIA
\item{[\ZIA]} C. Z\u{a}linescu, \slant Duality for vectorial nonconvex optimization by convexification and applications\endslant, An. St. Univ. Ia\c si, S. I Mat., {\bf 29}(1983), 15--34.   (This is available at $<$http://www.math.uaic.ro/~zalinesc/duality.pdf$>$.)  
\nmbr\ZZOR
\item{[\ZZOR]}-----, \slant Solvability results for sublinear functions and operators\endslant, Z. Oper. Res. Ser. A-B, {\bf 31}(1987), A79--A101. 
\nmbr\ZBOOK
\item{[\ZBOOK]}-----, \slant Convex analysis in
general vector spaces\endslant, (2002), World Scientific.
\Signoff
\bye